\documentclass{article}
\usepackage{amsfonts}
\usepackage{amsmath}
\usepackage{geometry}
\usepackage{amssymb}
\usepackage{graphicx}

\providecommand{\U}[1]{\protect\rule{.1in}{.1in}}
\providecommand{\U}[1]{\protect\rule{.1in}{.1in}}
\providecommand{\U}[1]{\protect\rule{.1in}{.1in}}

\begin{document}

\begin{center}
\bigskip

{\LARGE Superintegrable Cases of Four Dimensional Dynamical Systems}

\bigskip

O\u{g}ul Esen\footnote{%
oesen@yeditepe.edu.tr}$^{,a}$, Anindya Ghose Choudhury\footnote{%
aghosechoudhury@gmail.com}$^{,b}$, Partha Guha\footnote{%
partha@bose.res.in}$^{,c}$, Hasan G\"{u}mral\footnote{%
h.gumral@ack.edu.kw}\footnote{%
On leave of absence from Department of Mathematics, Yeditepe University,
Istanbul.}$^{,d}${\LARGE \ }

\bigskip
\end{center}

$^{a}$Department of Mathematics, Yeditepe University, 34755 Ata\c{s}ehir,
Istanbul, Turkey,

$^{b}$Department of Physics, Surendranath College, 24/2 Mahatma Gandhi Road,
Calcutta 700009, India,

$^{c}$SN Bose National Centre for Basic Sciences, JD Block, Sector III, Salt
Lake, Kolkata 700098, India,

$^{d}$Australian College of Kuwait, West Mishref, Kuwait

\bigskip

\bigskip

\textbf{Abstract:} Degenerate tri-Hamiltonian structures of Shivamoggi and
generalized Raychaudhuri equations are exhibited. For certain specific
values of the parameters, it is shown that hyperchaotic L\"{u} and Qi
systems are superintegrable and admit tri-Hamiltonian structures.

\bigskip

\textbf{Mathematics Classification (2010) }34C14, 34C20.

\smallskip


\bigskip

\section{Introduction}

Higher dimensional differential systems provide a setting and scope for
studies of diverse physical phenomena together with their associated
mathematical properties. In the case of integrable systems, for even
dimensions, one can investigate the existence of symplectic structures;
often in case of odd dimensions there exists what are known as Poisson
structures displaying a range of mathematical possibilities emerging from
the satisfaction of the Jacobi identity. Again from the dynamical point of
view higher dimensional systems form a basis for occurrence of chaos which
has become an intense field of study for a number of years now. Many
applications of such higher dimensional differential systems have been found
in the fields of engineering, biology, finance, medicine, climate studies
etc. A major hurdle in the analysis of higher-order ODEs arises from the
limited number of techniques at our disposal for their analysis. A majority
of such systems, which have applications in real life situations, are
nonlinear in character and until the advent of computational tools were
beyond the realm of concrete analysis by conventional methods. This,
however, should not be construed to undermine the power and utility of
analytical investigations of such systems, since they often provide a basis
for launching more focussed numerical investigations specially when dealing
with practical problems.

The analysis of higher-order ODEs is greatly aided by the existence of first
integrals since these allow for a reduction of the order of the system
through elimination of one or more variables. The determination of first
integrals is in itself a difficult task requiring, at times, substantial
guile as well as some amount of luck often in equal measure. While for
planar systems of ODEs there are some semi algorithmic methods depending on
the nature of the functions that are involved, the same cannot be said for
non-planar systems, thereby further complicating their analysis.
Nevertheless one can often adapt certain tools and techniques tailored for
planar systems and employ them for the analysis of non-planar systems.

The first nontrivial dimension for both the occurrence of chaos and non
canonical Hamiltonian formalism is the three dimensions where, the foremost
example of a chaotic systems, namely, the Lorenz system \cite{lorenz}, has
been observed to admit non canonical bi-Hamiltonian structure for certain
values of parameters involved, see \cite{GuNu93}. This structure has been
extended to integrable cases of Lotka-Volterra, May-Leonard, L\"{u} system,
Chen system etc. as well as to the Kermack-McKendrick model and some of its
generalizations describing the spread of epidemics \cite%
{ChLuLuYu,CGK,GuCh14,GuNu93,Nu90}. Important properties of three dimensional
Poisson structures surrounding all these globally integrable bi-Hamiltonian
systems are that the Jacobi identity is a scalar equation and is invariant
under multiplication with arbitrary function.

\subsection{Content of the work}

The purpose of the present work is to extend the integrability analysis to
four dimensions and to investigate the behavior of dynamical systems which
are known to be chaotic for generic values of their physical parameters.
This is essentially a non trivial problem with difficulties arising both
from integrability analysis and chaotic behavior. On one side, we have the
Jacobi identity which does not yield a single scalar equation, unlike the
corresponding situation for three-dimensional systems \cite{GC,Gum1,GuNu93}.
On the other side, the characterization of chaotic behavior of systems may
be more complicated. In four dimensions, one may encounter a chaotic system
with more than one positive Lyapunov exponent. These are usually defined as
hyperchaotic systems.

This work consists of two main sections. In the following section, we shall
provide a brief summary of the notion of integrability as enunciated by
Darboux. We shall recall method of Jacobi's last multiplier for planar
systems. Poisson structures and superintegrability will be reviewed.
Bi-Hamiltonian structures of superintegrable three dimensional dynamical
systems will be presented \cite{GuNu93}. As an example, Lorentz model will
be studied. Lorentz model has motivational importance for this study since
it is chaotic but, for some certain values of the parameters, it is
bi-Hamiltonian. This section will be ended by presenting a construction
which enables to write tri-Hamiltonian structure of a four dimensional
superintegrable dynamical system \cite{GoNu01}.

In section $3$, we shall exhibit superintegrable structures and
tri-Hamiltonian formulations of four dynamical systems in dimension four.
The first dynamical system is Shivamoggi equations governing a 4D
magneto-hydrodynamical system \cite{Shiva}. The second one is generalized
Raychaudhuri equations defining geodesic flows on surface of a deformable
media \cite{CGK,Ra55}. Both of the Shivamoggi and generalized Raychaudhuri
equations are autonomous and have time-independent first integrals.
Existence of three time-independent integrals enables us to construct
tri-Hamiltonian structures of these systems by simply following the general
procedure presented in subsection (\ref{4DP}). The third example is L\"{u}
system \cite{ChLuLuYu} and the fourth one is Qi system \cite{Qi}. They are
autonomous and hyperchaotic. We shall start with two time-dependent first
integrals for these systems. After a change of variables, the systems will
be nonautonumous whereas the integrals become time-independent. In the
intersection of level surfaces of two time-independent first integrals,
there will be a reduced autonomous dynamics in two dimensions. We shall
apply method of Jacobi's last multiplier for planar systems (c.f. subsection
(\ref{JLM})) to obtain a Hamiltonian function of the reduced dynamics. The
Hamiltonian functions for both of the reduced (L\"{u} and Qi) dynamics turn
out to be time-dependent. We shall change the time parameter and put some
conditions on the other parameters in order to obtain an additional
time-independent first integral. Obtaining three time-independent first
integrals, we shall present tri-Hamiltonian structures of L\"{u} and Qi
systems. Thus we shall be able to show that even though L\"{u} and Qi
systems are hyperchaotic, for certain values of the parameters, they are
superintegrable and admit tri-Hamiltonian structures.

\section{Integrals, Superintegrability and Multi-Hamiltonian Systems}

\subsection{The first and second integrals}

A first integral of a system of ODE's%
\begin{equation}
\dot{x}^{a}=X^{a}(t,x^{0},x^{1},...,x^{n-1})\text{, \ \ \ }a=0,1,...,n-1,
\label{eom}
\end{equation}
is a differentiable function $I(t,x^{0},x^{1},...,x^{n-1})$ that retains a
constant value on any integral curve of the system, that is, its derivative
with respect to time $t$ vanishes on the solution curves of the system (\ref%
{eom}), see \cite{Go01,OLV}. After introducing a vector field
\begin{equation}
X(t,x^{0},x^{1},...,x^{n-1})=X^{a}(t,x^{0},x^{1},...,x^{n-1})%
\partial_{x^{a}}+\partial_{t},
\end{equation}
a first integral can be defined by the condition $X(I)=0$. When the
coefficient functions $X^{a}$ in (\ref{eom}) have no explicit dependence on
the time $t$, the system is called autonomous. In this case, a first
integral may still contain $t$ explicitly.

A second integral is a $C^{1}$ function $J$ satisfying the identity
\begin{equation*}
X(J)=\vartheta J.
\end{equation*}
Here, $\vartheta$ is a real valued function called cofactor of the second
integral. Polynomial second integrals for (polynomial) vector fields, called
Darboux polynomials, simplify the determination of the first integrals \cite%
{Da878}. When we have two relatively prime Darboux polynomials $P_{1}$ and $%
P_{2}$, having a common cofactor, of a (polynomial) vector field $X$, then
the fraction $P_{1}/P_{2}$ is a rational first integral of $X$. The inverse
of this statement is also true, that is if we have a rational first integral
$P_{1}/P_{2}$ of a vector field $X$, then $P_{1}$ and $P_{2}$ are Darboux
polynomials for $X$. For planar (polynomial) vector fields, when we have a
certain number of relatively prime irreducible Darboux polynomials, not
necessarily having a common cofactor, it is possible to write first
integrals from the Darboux polynomials \cite{Da878,DuLlAr06,Ju79,Si92}. The
generalization of this theorem for arbitrary dimensions is still open.

In \cite{PrSi83,Si92}, a semi-algorithm, called Prelle-Singer method, is
presented for the determinations of elementary first integrals for planar
systems. For higher dimensional systems, Darboux polynomials are useful,
though sometimes the use of a specific ansatz or a polynomial in one
variable (of particular degree) with coefficients depending on the remaining
variables remains the only option. One may also use a variant of the
Prelle-Singer/Darboux method to derive what are called quasi-rational first
integrals \cite{Ma94}.

\subsection{Method of Jacobi's last multiplier \label{JLM}}

The Jacobi's last multiplier is a useful tool for deriving an additional
first integral for a system of $n$ first-order ODEs when $n-2$ first
integrals of the system are known. It allows us to determine the
Lagrangian/Hamiltonian of a planar system in many cases \cite%
{Go01,Jac1,Jac2,Whi}. In recent years a number of articles have dealt with
this particular aspect \cite{CGK,Nucci,Leach}.

For the case $n=2$, Jacobi's last multiplier lets us to find a
transformation $(x,y,t)\rightarrow(Q,P,t)$ such that the system
\begin{equation}
\dot{x}=f(x,y,t),\;\;\;\dot{y}=g(x,y,t).  \label{5.4}
\end{equation}
may be expressed in the form of Hamilton's equations%
\begin{equation}
\dot{Q}=\partial_{P}H,\;\;\;\dot{P}=-\partial_{Q}H.  \label{5.6}
\end{equation}
The method goes as follows. Wedge product of the Cartan one-forms $\theta
^{x}=dx-f(x,y,t)dt\;$and$\;\theta^{y}=dy-g(x,y,t)dt$ is
\begin{equation}
\omega=\theta^{x}\wedge\theta^{y}.  \label{alpha}
\end{equation}
The Hamiltonian formulation in the extended phase space $%
\mathbb{R}
^{2}\times%
\mathbb{R}
^{2}$ is characterized by the Poincar\'{e}-Cartan one-form $\Theta=PdQ-Hdt$
from which we may define a closed two-form $\Omega$ by minus the exterior
derivative of $\Theta$, that is,
\begin{equation}
\Omega=dQ\wedge dP+dH\wedge dt.  \label{5.8}
\end{equation}
It follows that, the two-form $\omega$ in Eq.(\ref{alpha})\ is proportional
to $\Omega$ with
\begin{equation}
M\omega=\Omega  \label{alpha-omega}
\end{equation}
for some function $M$, called Jacobi Last Multiplier. Since $\Omega$ is
closed, $M\omega$ must be closed. This leads to the following differential
equation
\begin{equation}
\partial_{t}M+\partial_{x}(Mf)+\partial_{y}(Mg)=0  \label{5.10}
\end{equation}
determining $M$. After solving equation (\ref{5.10}), the canonical
coordinates can be obtained by substituting $M$ into Eq.(\ref{alpha-omega}).
To determine the Hamiltonian function $H$, there are two possible cases
depending whether $\partial_{t}M$ vanishes or not.

If the multiplier $M$ does not depend on the time $t$ then the first term on
the left hand side of Eq.(\ref{5.10}) drops. In this case, from the equality
$M\omega=\Omega$, we arrive at%
\begin{equation}
M\left( fdy-gdx\right) =dH,  \label{5.11}
\end{equation}
which relates $M$ with the Hamiltonian function $H$.

When the multiplier $M$ depends on time explicitly, we introduce two
auxiliary functions $\phi$ and $\psi$ such that the structure of Eq.(\ref%
{5.11}) is preserved, that is, for a pair of functions $\psi$ and $\phi$ we
have
\begin{equation}
M((f-\psi)dy-(g-\phi)dx)=dH+\vartheta dt,  \label{5.13}
\end{equation}
for a real valued function $\vartheta$ of $(Q,P,t)$, \cite{Ca09}. This
occurs if the condition
\begin{equation}
\partial_{x}(M(f-\psi))+\partial_{y}(M(g-\phi))=0  \label{5.12}
\end{equation}
is satisfied. By adding and subtracting two-form $M\psi dy\wedge dt-M\phi
dx\wedge dt$ into Eq.(\ref{alpha-omega}), we obtain%
\begin{align*}
M\omega & =M\theta^{x}\wedge\theta^{y}=Mdx\wedge dy+M(fdy-gdx)\wedge dt \\
& =Mdx\wedge dy+M(fdy-gdx)\wedge dt\pm\left( M\psi dy\wedge dt-M\phi
dx\wedge dt\right) \\
& =M(dx-\psi dt)\wedge(dy-\phi dt)+M\left( (f-\psi)dy-(g-\phi)dx\right)
\wedge dt,
\end{align*}
where the first two-form in the last line is the symplectic two-form
\begin{equation}
M(dx-\psi dt)\wedge(dy-\phi dt)=dQ\wedge dP,  \label{Canoor}
\end{equation}
and the latter can be obtained by taking the exterior product of both sides
of Eq.(\ref{5.13}) by $dt.$ Note that, substitutions of the auxiliary
functions $\psi$ and $\phi$, and the multiplier $M$ into Eq.(\ref{5.13})
will enable us to find the Hamiltonian function whereas substitutions into
Eq.(\ref{Canoor}) will determine the canonical coordinates.

\subsection{Poisson structures and superintegrability}

A Poisson structure on an $n$-dimensional space is a skew-symmetric bracket $%
\left\{ \text{ , }\right\} $ on the space of real-valued smooth functions
satisfying the Leibnitz and the Jacobi identities \cite{LM, LaPi12,OLV,
wei83}. We define the Poisson bracket of two functions $F$ and $H$ by
\begin{equation}
\left\{ F,H\right\} =\nabla F\cdot N\nabla H,  \label{PB}
\end{equation}
where $\nabla F$ and $\nabla H$ are gradients of the functions $F$ and $H$
respectively, and $N$ is Poisson matrix. The Poisson bracket (\ref{PB})
automatically satisfies the Leibnitz identity. In a local chart ($x^{a}$),
the Jacobi identity takes the form
\begin{equation}
N^{a[b}\partial_{x^{a}}N^{cd]}=0,  \label{jcb}
\end{equation}
where $N^{ab}$ are components of the Poisson matrix $N$ and, $[$ $\ ]$
refers anti-symmetrization. The Jacobi identity (\ref{jcb}) is trivially
satisfied in two dimensions, it gives a scalar equation in three dimensions
and, it results in four equations in four dimensions.

A dynamical system is Hamiltonian if it can be written as
\begin{equation}
\dot{x}=N\nabla H  \label{HamEq}
\end{equation}%
for $H$ being a real valued function, called Hamiltonian function, and $N$
being a Poisson matrix. For autonomous Hamiltonian systems, the Hamiltonian
function is conserved, that is it is also a first integral of the system. It
is also possible to find a Hamiltonain formulation of a nonautonomous
system, but in this case the Hamiltonian is not a constant of motion, hence
not a first integral of the system \cite{AbMa78}.

A Hamiltonian system in $n$ dimensions is called maximally superintegrable
if there are $n-1$ first integrals. Existence of $n-1$ integrals lets one to
reduce the systems of differential equations to one quadrature. For the case
of three dimensions, two first integrals are required for maximal
superintegrability whereas for four dimensions, three first integrals are
needed.

\subsection{Superintegrability, bi-Hamiltonian systems and Lorentz example}

A dynamical system is bi-Hamiltonian if it admits two different Hamiltonian
structures
\begin{equation*}
\dot{x}=N^{\left( 1\right) }\nabla H_{2}=N^{\left( 2\right) }\nabla H_{1}
\end{equation*}%
such that any linear pencil $N^{\left( 1\right) }+cN^{\left( 2\right) }$ of
Poisson matrices satisfies the Jacobi identity \cite%
{AbGu09,BlWo89,MaMo84,OLV}. In this case, one can generate recursively
enough constants of motion to ensure integrability \cite{Fe94,Ol90}. In \cite%
{GuNu93}, it is shown that, Hamiltonian structures of dynamical systems in
three dimensions always come in compatible pairs to form a bi-Hamiltonian
structure with Poisson matrices
\begin{equation*}
N^{(i)ab}=-\epsilon ^{ij}\epsilon ^{abc}\partial _{x^{c}}H_{j}\text{, \ \ }%
a,b,c=1,2,3,\text{ \ \ }i,j=1,2,
\end{equation*}%
where $\epsilon ^{ij}$ and $\epsilon ^{abc}$ are completely antisymmetric
tensors of rank two and three, respectively. In the rest of this subsection,
we, particularly, focus on Lorentz system. Although it is chaotic, for some
certain values of its parameters, it is superintegrable and bi-Hamiltonian
\cite{GuNu93}. In that sense, Lorentz system has a motivational importance
for this paper.

The Lorenz model \cite{lorenz} is a three component dynamical system%
\begin{equation}
\dot{x}=\sigma(y-x)\,,\text{ }\dot{y}\text{\ }=\rho x-xz-y,\text{ \ \ \ }%
\dot{z}=-\beta z+xy  \label{lor}
\end{equation}
where $\sigma$ and $\rho$ are the Prandtl and Rayleigh numbers,
respectively, and $\beta$ is another dimensionless number, the aspect ratio.
It exhibits chaotic behavior for most values of these parameters. However,
it admits bi-Hamiltonian structure in two limits which can, most
conveniently, be characterized by the Rayleigh number $\rho$, namely $\rho=0$
and $\rho\rightarrow\infty$.

The case $\rho=0,$ $\sigma=1/2,$ $\beta=1$ is known to admit two
time-dependent conserved quantities \cite{segur,tw,kus,st}. The
transformation
\begin{equation}
x=\frac{1}{2}{\bar{t}}u\,,\text{ \ \ }y=\frac{1}{4}{\bar{t}}^{2}v,\text{ \ \
}z=\frac{1}{4}{\bar{t}}^{2}w,\text{ \ \ }t=-\log({\bar{t}}^{2}/4)
\end{equation}
of dynamical variables and time brings the Lorenz system to the form
\begin{equation}
u^{\prime}=\frac{1}{2}v\,,\text{ \ \ }v^{\prime}=-uw,\text{ \ \ }%
w^{\prime}=uv  \label{has}
\end{equation}
with prime denoting differentiation with respect to ${\bar{t}}$. In terms of
the new dynamical variables we find that%
\begin{equation}
H_{1}=w-u^{2}\text{, \ \ \ }H_{2}=v^{2}+w^{2},  \label{tindepc}
\end{equation}
are conserved. The Hamiltonian structure functions are given by
\begin{equation}
N^{(1)}=\frac{1}{4}\left(
\begin{array}{ccc}
0 & 1 & 0 \\
-1 & 0 & -2u \\
0 & 2u & 0%
\end{array}
\right) ,\text{ \ \ \ }N^{(2)}=\frac{1}{2}\left(
\begin{array}{ccc}
0 & -w & v \\
w & 0 & 0 \\
-v & 0 & 0%
\end{array}
\right) .  \label{j1}
\end{equation}

The second completely integrable case of the Lorenz system is the
conservative limit \cite{tw},\cite{se} obtained through the scaling of Eq.(%
\ref{lor}) with%
\begin{equation}
t\rightarrow\epsilon t,\text{ \ \ }x\rightarrow\frac{1}{\epsilon}x,\text{ \
\ }y\rightarrow\frac{1}{\sigma{\epsilon}^{2}}y,\text{ \ \ }z\rightarrow
\frac{1}{\sigma{\epsilon}^{2}}z,\text{ \ \ }\epsilon=\frac{1}{(\sigma
\rho)^{1/2}}  \label{trclim}
\end{equation}
and taking the limit $\epsilon\rightarrow0.$ This brings Eq.(\ref{lor}) to
the form%
\begin{equation}
\dot{x}=y\,,\text{ \ }\dot{y}\text{\ }=-xz+x,\text{ \ \ \ }\dot{z}=xy
\label{clor}
\end{equation}
possessing two first integrals%
\begin{equation}
H_{1}=\frac{1}{2}(y^{2}+z^{2}-x^{2}),\text{ \ \ \ }H_{2}=\frac{1}{2}x^{2}-z.
\label{hc12}
\end{equation}
The Poisson structures are defined by the matrices
\begin{equation}
N^{\left( 1\right) }=\left(
\begin{array}{ccc}
0 & z & -y \\
-z & 0 & -x \\
y & x & 0%
\end{array}
\right) ,\text{ \ \ \ }N^{\left( 2\right) }=\left(
\begin{array}{ccc}
0 & 1 & 0 \\
-1 & 0 & -x \\
0 & x & 0%
\end{array}
\right)  \label{j1l}
\end{equation}
which are compatible.

We have shown that, in two extreme cases Lorentz system is completely
integrable and admits bi-Hamiltonian structure. In \cite{GuNu93}, it was
shown that, these two cases are the same, hence they can be recognized as an
example of a Nambu system \cite{nam} because the equations of motion take
the form
\begin{equation}
\dot{x}^{i}=\epsilon^{ijk}\partial_{j}H_{1}\partial_{k}H_{2}  \label{nambu}
\end{equation}
where $\epsilon^{ijk}$ is the completely skew Levi-Civita tensor.

\subsection{Superintegrability and tri-Hamiltonian systems in four
dimensions \label{4DP}}

In \cite{GoNu01}, a construction of multi-Hamiltonian formalism of maximal
superintegrable systems is presented. For the case of four dimensions, this
construction results with a tri-Hamiltonian system%
\begin{equation}
\dot{x}=\vartheta N^{\left( 1\right) }\bar{\nabla}H_{1}=\vartheta N^{\left(
2\right) }\bar{\nabla}H_{2}=\vartheta N^{\left( 3\right) }\bar{\nabla}H_{3}
\label{tri-Ham}
\end{equation}
consisting of three compatible Poisson matrices $N_{1}$, $N_{2}$ and $N_{3}$
with linearly independent Hamiltonians $H_{1},$ $H_{2}$ and $H_{3}$,
respectively. Here, $\bar{\nabla}$ denotes the gradient vector in four
dimension. The real valued function $\vartheta$ in Eq.(\ref{tri-Ham}) is
called the conformal factor. In general, by multiplying a Poisson matrix $N$
with a real valued function $\vartheta$, we obtain a new Poisson matrix $%
\vartheta N$ in two and three dimensions. In dimension four, $\vartheta N$
is a Poisson matrix if $N$ is degenerate, that is the rank of $N$ fails to
be full. We shall refer this as conformal invariance.

A Poisson matrix can be written as%
\begin{equation*}
N=%
\begin{pmatrix}
0 & -U^{1} & -U^{2} & -U^{3} \\
U^{1} & 0 & -V^{3} & V^{2} \\
U^{2} & V^{3} & 0 & -V^{1} \\
U^{3} & -V^{2} & V^{1} & 0%
\end{pmatrix}%
\end{equation*}
which enables us to define three-component vector functions $\mathbf{U}%
=\left( U^{1},U^{2},U^{3}\right) $ and\ $\mathbf{V}=\left(
V^{1},V^{2},V^{3}\right) $ of four variables $(u,x,y,z)=\left( u,\mathbf{x}%
\right) .$ Accordingly, we denote the gradient vector as follows%
\begin{equation*}
\bar{\nabla}=\left( \partial_{u},\nabla\right) =\left( \partial
_{u},\partial_{x},\partial_{y},\partial_{z}\right) ,
\end{equation*}
where $\nabla=\left( \partial_{x},\partial_{y},\partial_{z}\right) $. The
Jacobi identity in Eq.(\ref{jcb}) gives four equations that can be divided
into one scalar equation and one vector equation%
\begin{align}
\partial_{u}(\mathbf{U}\cdot\mathbf{V}) & =\mathbf{V}\cdot\left( \partial_{u}%
\mathbf{U}-\nabla\times\mathbf{V}\right) ,  \label{dsjk} \\[0.08in]
\nabla(\mathbf{U}\cdot\mathbf{V}) & =\mathbf{V}(\nabla\cdot\mathbf{U})-%
\mathbf{U}\times\left( \partial_{u}\mathbf{U}-\nabla\times\mathbf{V}\right) .
\label{dvjk}
\end{align}
Note that, the left-hand sides of Eqs.(\ref{dsjk}) and (\ref{dvjk}) vanish
for degenerate $(\mathbf{U}\cdot\mathbf{V}=0)$ matrices.

Let $H_{1}$ and $H_{2}$ be two time-independent first integrals, we define
three-component vector functions%
\begin{equation}
\mathbf{U}=\nabla H_{1}\times\nabla H_{2}\text{, \ \ \ }\mathbf{V}%
=\partial_{u}H_{1}\nabla H_{2}-\partial_{u}H_{2}\nabla H_{1}.  \label{UandV}
\end{equation}
Note that, the vector functions $\mathbf{U}$ and $\mathbf{V}$ in Eq.(\ref%
{UandV}) are orthogonal which implies $N$ is degenerate for arbitrary $C^{1}$
functions $H_{1}$ and $H_{2}$. For Jacobi identity expressed in Eqs.(\ref%
{dsjk}) and (\ref{dvjk}), we compute%
\begin{align}
\nabla\cdot\mathbf{U} & =\nabla\cdot\left( \nabla H_{1}\times\nabla
H_{2}\right) =\nabla\cdot\left( \nabla\times H_{1}\nabla H_{2}\right) =0
\notag \\
\partial_{u}\mathbf{U}-\nabla\times\mathbf{V} & =\partial_{u}\left( \nabla
H_{1}\times\nabla H_{2}\right) -\nabla\times\left( \partial_{u}H_{1}\nabla
H_{2}-\partial_{u}H_{2}\nabla H_{1}\right)  \label{cak} \\
& =\nabla\left( \partial_{u}H_{1}\right) \times\nabla H_{2}+\nabla
H_{1}\times\nabla\left( \partial_{u}H_{2}\right) -\nabla\left( \partial
_{u}H_{1}\right) \times\nabla H_{2}+\nabla\left( \partial_{u}H_{2}\right)
\times\nabla H_{1}=\mathbf{0}.  \notag
\end{align}
For a given four component dynamical vector field $X=(X^{u},\mathbf{X})$,
the Hamilton's equations (\ref{HamEq}) take the particular form
\begin{align*}
\dot{u} & =X^{u}=-\left( \nabla H_{1}\times\nabla H_{2}\right) \cdot\nabla H,
\\
\mathbf{\dot{x}} & =\mathbf{X=(}\nabla H_{1}\times\nabla H_{2})\partial
_{u}H+(\nabla H_{2}\times\nabla H)\partial_{u}H_{1}+(\nabla H\times\nabla
H_{1})\partial_{u}H_{2}.
\end{align*}
One can obtain the conservation equations%
\begin{equation}
X^{u}\partial_{u}H+\mathbf{X}\cdot\nabla H=X^{u}\partial_{u}H_{1}+\mathbf{X}%
\cdot\nabla H_{1}=X^{u}\partial_{u}H_{2}+\mathbf{X}\cdot\nabla H_{2}=0
\label{ConsLaw}
\end{equation}
for Hamiltonian function $H$ as well as for functions $H_{1}$ and $H_{2}$
defining the Poisson structure. The conservation laws (\ref{ConsLaw}) for $%
H_{1}$ and $H_{2}$ show that these functions are Casimirs of the Poisson
matrix $N$ they constructed, that is $N\bar{\nabla}H_{i}=\mathbf{0}$ for $%
i=1,2$.

By interchanging cyclically the roles of the functions $\left(
H_{1},H_{2},H=H_{3}\right) $ one obtains tri-Hamiltonian structure of the
system. In a compact notation, the coefficients of Poisson matrices are in
form
\begin{equation}
N^{(i)ab}=-\epsilon^{ijk}\epsilon^{abcd}\partial_{x^{c}}H_{j}%
\partial_{x^{d}}H_{k}\text{, \ }a,b,c,d=0,1,2,3,\text{ \ }i,j,k=1,2,3,
\label{supern}
\end{equation}
where $\epsilon^{ijk}$ and $\epsilon^{abcd}$ are completely antisymmetric
tensors of rank three and four, respectively. In the vector notation, let $%
\left( \mathbf{U}_{1},\mathbf{V}_{1}\right) ,\left( \mathbf{U}_{2},\mathbf{V}%
_{2}\right) $ and $\left( \mathbf{U}_{3},\mathbf{V}_{3}\right) $ be
three-component vector functions of three Poisson matrices $N^{\left(
1\right) }$, $N^{\left( 2\right) }$ and $N^{\left( 3\right) }$ given by
\begin{align}
\mathbf{U}_{1} & =\nabla H_{2}\times\nabla H_{3},\text{ \ \ }\mathbf{V}%
_{1}=\partial_{u}H_{2}\nabla H_{3}-\partial_{u}H_{3}\nabla H_{2},  \label{N}
\\
\mathbf{U}_{2} & =\nabla H_{3}\times\nabla H_{1},\text{ \ \ }\mathbf{V}%
_{2}=\partial_{u}H_{3}\nabla H_{1}-\partial_{u}H_{1}\nabla H_{3},  \notag \\
\mathbf{U}_{3} & =\nabla H_{1}\times\nabla H_{2},\text{ \ \ }\mathbf{V}%
_{3}=\partial_{u}H_{1}\nabla H_{2}-\partial_{u}H_{2}\nabla H_{1},  \notag
\end{align}
respectively. It is straight forward to verify that the conditions%
\begin{equation}
\nabla\cdot\mathbf{U}_{i}=0\text{, \ \ \ \ \ }\frac{\partial\mathbf{U}_{i}}{%
\partial u}-\nabla\times\mathbf{V}_{i}=\mathbf{0}  \label{jac}
\end{equation}
are satisfied to guarantee Jacobi identities for all Poisson matrices $%
N^{\left( i\right) }$, $i=1,2,3$. It is also straightforward to see that%
\begin{equation}
\Lambda_{ij}=\mathbf{U}_{i}\cdot\mathbf{V}_{j}+\mathbf{U}_{j}\cdot \mathbf{V}%
_{i}=0  \label{lam}
\end{equation}
which shows that all three Hamiltonian structures are mutually compatible.
Thus, Poisson structures for superintegrable systems in dimension four
always form compatible pairs.

\section{Examples}

\subsection{Shivamoggi equations}

Shivamoggi equations are arising in the context of four dimensinal
magnetohydrodynamics and are given by%
\begin{equation}
\dot{u}=-uy,\;\;\;\dot{x}=zy,\;\ \;\dot{y}=zx-u^{2},\;\;\;\dot{z}=xy,
\label{SE}
\end{equation}
see \cite{GuCh14,Shiva}. The first integrals of this system of equations are
\begin{equation}
H_{1}=x^{2}-z^{2},\;\;\;H_{2}=z^{2}+u^{2}-y^{2},\;\;\;H_{3}=u(z+x).
\label{FISE}
\end{equation}
From Eq.(\ref{N}), we identify the vectors $\mathbf{U}_{i}$ and $\mathbf{V}%
_{i}$ of Poisson matrices $N^{\left( i\right) }$ for the Hamiltonian
functions $H_{1},H_{2}$\ and $H_{3}$, $i=1,2,3$, we find%
\begin{align}
\mathbf{U}_{1} & =2u\left( -y,z,y\right) ,\text{ \ \ \ }\mathbf{V}%
_{1}=2\left( u^{2},y(x+z),u^{2}-z(x+z)\right) ,  \notag \\
\mathbf{U}_{2} & =2\left( x+z\right) \left( 0,u,0\right) ,\text{ \ \ \ }%
\mathbf{V}_{2}=\left( x,0,-z\right) ,  \notag \\
\mathbf{U}_{3} & =-4\left( yz,zx,xy\right) ,\text{ \ \ \ }\mathbf{V}%
_{3}=4u\left( x,0,-z\right) ,  \label{PoiShi}
\end{align}
respectively. Note that, all of these three Poisson matrices are degenerate,
since $\mathbf{U}_{i}\cdot\mathbf{V}_{i}=0$ holds for all $i=1,2,3$. The
equations of motion can be written as%
\begin{equation*}
X=\vartheta N^{\left( 1\right) }\bar{\nabla}H_{1}=\vartheta N^{\left(
2\right) }\bar{\nabla}H_{2}=\vartheta N^{\left( 3\right) }\bar{\nabla}H_{3}%
\text{, \ \ \ }\vartheta=-\frac{1}{4(x+z)}
\end{equation*}
up to multiplication with a conformal factor $\vartheta$ for all three.

The quadratic system (\ref{SE}) may be related to certain four dimensional
Lie algebra by defining a linear Poisson structure \cite{LaPi12}. Using the
first two first integrals, it is possible to define the Hamiltonian
\begin{equation}
H=H_{1}-H_{2}=x^{2}+y^{2}-2z^{2}-u^{2},
\end{equation}
for the system (\ref{SE}) satisfying $\dot{X}=N\bar{\nabla}H$. Here, $N$ is
a Poisson matrix defined by $\mathbf{U}=(0,y,0)$ while $\mathbf{V}%
=(-x,0,-2u) $. Note that, $\mathbf{U}\cdot\mathbf{V}=0$ so that the this
Poisson structure is also degenerate.

\subsection{ Generalized Raychaudhuri equations}

The generalized Raychaudhuri equations are a set of coupled first-order ODEs
related to geodesic flows on surface of a deformable media \cite{Ra55}. In
the case of a two dimensional curved surface of constant curvature they give
rise to the following set of equations
\begin{align}
\dot{x}+\frac{1}{2}x^{2}+\alpha x+2(y^{2}+z^{2}-u^{2})+2\beta & =0,  \notag
\\
\dot{y}+(\alpha +x)y+\gamma & =0,  \notag \\
\dot{z}+(\alpha +x)z+\delta & =0,  \notag \\
\dot{u}+(\alpha +x)u& =0
\end{align}%
upon use of the exact solutions of geodesic equations. The vector field%
\begin{equation}
X=-\left( \frac{1}{2}x^{2}+2(y^{2}+z^{2}-u^{2})\right) \partial
_{x}-xy\partial _{y}-xz\partial _{z}-xu\partial _{u}  \label{RE}
\end{equation}%
is associated to a particular case obtained by setting all the four
parameters $\alpha =\beta =\gamma =\delta =0$. This system admits the
following Darboux polynomials
\begin{equation*}
J_{1}=y,\text{ \ }J_{2}=z,\text{ \ }J_{3}=u,\text{ \ }%
J_{4}=y^{2}+z^{2}-u^{2}-\frac{1}{4}x^{2}
\end{equation*}%
together with cofactors
\begin{equation*}
\lambda _{1}=\lambda _{2}=\lambda _{3}=\lambda _{4}=-x,
\end{equation*}%
respectively, see \cite{CGK,Valls}. This is the case where the Darboux
polynomails have common cofactor $-x$, hence we have several first integrals
obtained by the fractions of the Darboux polynomials. For example,
\begin{equation}
H_{1}=\frac{z}{u}\text{, \ \ }H_{2}=\frac{y}{z},\text{ \ \ }H_{3}=\frac{1}{u}%
\left( y^{2}+z^{2}-u^{2}-\frac{1}{4}x^{2}\right) ,\ \ H_{4}=\frac{1}{y}%
\left( y^{2}+z^{2}-u^{2}-\frac{1}{4}x^{2}\right) .  \label{FIRE}
\end{equation}

The three-component vector functions defined in Eqs.(\ref{N}) for Poisson
matrices corresponding to the Hamiltonian functions $H_{1},$ $H_{2},$ $H_{3}$
are
\begin{align*}
\mathbf{U}_{1} & =\frac{-2}{uz^{2}}\left( 4\left( y^{2}+z^{2}\right)
,xy,xz\right) ,\text{ \ \ \ }\mathbf{V}_{1}=\frac{4}{z^{2}}\left(
y^{2}+z^{2}+u^{2}-\frac{1}{4}x^{2}\right) \left( 0,-z,y\right) , \\
\mathbf{U}_{2} & =\frac{-2}{u^{2}}\left( 4y,x,0\right) ,\text{ \ \ \ }%
\mathbf{V}_{2}=\frac{2}{u^{3}}\left( xz,-4yz,2u^{3}\left( y^{2}+z^{2}+u^{2}-%
\frac{1}{4}x^{2}\right) -4z^{2}\right) , \\
\mathbf{U}_{3} & =\frac{1}{uz}\left( -1,0,0\right) ,\text{ \ \ \ }\mathbf{V}%
_{3}=\frac{1}{zu^{2}}\left( 0,-z,y\right) ,
\end{align*}
respectively. The equations of motion can be written as%
\begin{equation*}
X=\vartheta N^{(1)}\bar{\nabla}H_{1}=\vartheta N^{\left( 2\right) }\bar{%
\nabla}H_{2}=\vartheta N^{\left( 3\right) }\bar{\nabla}H_{3}\text{, \ \ \ }%
\vartheta=-\frac{1}{2}zu^{3}
\end{equation*}
where the conformal factors $\vartheta$ are the same for all three.

It may be instructive to start with two Hamiltonian functions $H_{1}\ $and $%
H_{2}$ in Eq.(\ref{FIRE}), and obtain $H_{3}$ using the method of Jacobi's
last multiplier. On the common level sets of $H_{1}=\kappa\ $and $H_{2}=\tau$%
, the dynamics generated by the vector field (\ref{RE}) reduces to the
dynamics%
\begin{equation}
\dot{x}=-\frac{1}{2}x^{2}+\mu z^{2}\text{,\ \ \ }\dot{z}=-xz,  \label{RedRay}
\end{equation}
where we have eliminated $y=\tau z$ and $u=z/\kappa$ in favour of $x$ and $%
z, $ and $\mu=2/\kappa^{2}-2\tau^{2}-2$ is a constant for dynamics of (\ref%
{RedRay}). By solving the defining PDE in Eq.(\ref{5.10}) for the case of
the reduced dynamics in Eq.(\ref{RedRay}), we find the multiplier $M=1/z^{2}$%
. This enables us to determine the Hamiltonian function
\begin{equation}
H\left( x,z\right) =\frac{x^{2}}{2z}+\mu z.  \label{HonRay}
\end{equation}
The canonical coordinates are $Q=x$ and $P=-1/z$, hence the symplectic two
form is $\Omega=\left( 1/z^{2}\right) dx\wedge dz.$ Substitutions of $\mu,$ $%
\kappa$ and $\tau$ into the Hamiltonian function $H$ in Eq.(\ref{HonRay})
result a scalar multiple of $H_{3}$ in Eq.(\ref{FIRE}) as expected.

Corresponding Poisson matrix for the fourth integral $H_{4}$ in Eq.(\ref%
{FIRE}) may be in the form%
\begin{equation*}
N=lN^{\left( 1\right) }+mN^{\left( 2\right) }+nN^{\left( 3\right) }
\end{equation*}
for arbitrary functions $l,m$ and $n$. We find that these functions are
subjected to the condition%
\begin{equation*}
-\frac{8}{zu^{3}}(lH_{2}H_{3}+mH_{1}H_{3}+nH_{1}H_{2})=1.
\end{equation*}

\subsection{Hyperchaotic L\"{u} system \label{FIHLS}}

We consider a hyperchaotic system of four first-order ODEs obtained from the
L\"{u} system \cite{ChLuLuYu} by adding an additional variable. The set of
equations is%
\begin{align}
\dot{u} & =\delta u+xz,  \notag \\
\dot{x} & =\alpha(y-x)+u,  \notag \\
\dot{y} & =\gamma y-xz,  \notag \\
\dot{z} & =-\beta z+xy,  \label{LuG}
\end{align}
where $\alpha$, $\beta$, $\gamma$, and $\delta$ are real constant
parameters. When the parameters $\gamma=-\beta=\delta$, we obtain two
time-dependent first integrals
\begin{equation}
I_{1}=e^{-\gamma t}(y+u),\;\;\;I_{2}=e^{-2\gamma t}(y^{2}+z^{2}).
\label{ILu1}
\end{equation}
We change the variables $\left( u,x,y,z\right) $ with $\left( s,p,q,r\right)
$\ according to%
\begin{equation}
s=ue^{-\gamma t},\text{ \ \ }p=xe^{\alpha t},\;\;\;q=ye^{-\gamma
t},\;\;\;r=ze^{-\gamma t}.  \label{cov}
\end{equation}
In this new coordinates, the L\"{u} system becomes nonautonomous
\begin{equation}
\dot{s}=rpe^{-\alpha t},\;\;\;\dot{p}=(\alpha q+s)e^{(\alpha+\gamma )t},\;\;%
\dot{q}=-rpe^{-\alpha t},\;\;\;\dot{r}=qpe^{-\alpha t},  \label{LunonAut}
\end{equation}
whereas the first integrals $I_{1}$ and $I_{2}$ in Eq.(\ref{cov}) become
time-independent
\begin{equation}
H_{1}\left( s,p,q,r\right) =q+s,\;\;\;H_{2}\left( s,p,q,r\right)
=q^{2}+r^{2}.  \label{ILu}
\end{equation}

On the intersection of the level hypersurfaces $H_{1}=\kappa$ and $%
H_{2}=\tau $, the set of equations (\ref{LunonAut}) reduces to a
nonautonomous planar system%
\begin{equation}
\dot{p}=(\kappa+(\alpha-1)q)e^{(\alpha+\gamma)t},\;\;\;\dot{q}=-p\left(
\tau-q^{2}\right) ^{\frac{1}{2}}e^{-\alpha t},  \label{Lu3}
\end{equation}
where we have eliminated $s=(\kappa-q)$ and $r=\sqrt{\tau-q^{2}}$ in favor
of $p$ and $q$. It is straightforward to verify that the Jacobi's last
multiplier for the system (\ref{Lu3}) is
\begin{equation*}
M=\left( \tau-q^{2}\right) ^{-\frac{1}{2}}.
\end{equation*}
Consequently, we write a Hamiltonian function
\begin{equation}
H\left( q,p\right) =e^{(\alpha+\gamma)t}\left[ \kappa\arcsin\left( \frac{q}{%
\sqrt{\tau}}\right) -(\alpha-1)\sqrt{\tau-q^{2}}\right] +\frac {1}{2}%
p^{2}e^{-\alpha t}.  \label{HamLu1}
\end{equation}
In terms of the canonical variables%
\begin{equation}
Q=\arcsin\left( \frac{q}{\sqrt{\tau}}\right) ,\;\;\;P=p
\end{equation}
nonautonomous system (\ref{Lu3}) takes the form
\begin{equation}
\dot{Q}=Pe^{-\alpha t}\text{, \ \ }\dot{P}=-e^{(\alpha+\gamma)t}\left(
\kappa+(\alpha-1)\sqrt{\tau}\sin Q\right)  \label{Lu2DCan}
\end{equation}
whereas the Hamiltonian $H$ in Eq.(\ref{HamLu1}) becomes
\begin{equation}
H=e^{(\alpha+\gamma)t}\left( \kappa Q-(\alpha-1)\sqrt{\tau}\cos Q\right) +%
\frac{1}{2}P^{2}e^{-\alpha t},  \label{HamLu2}
\end{equation}
so that, we write the system (\ref{Lu2DCan}) in form of Hamilton's equations
$\dot{Q}=\partial H/\partial P$ and $\dot{P}=\partial H/\partial Q$. The
system (\ref{Lu2DCan}) is nonautonomuous, hence the time-dependent
Hamiltonian $H$ in (\ref{HamLu2})\ is not an integral invariant of the
motion because its total derivative $dH/dt$ with respect to time $t$ equals
to $\partial H/\partial t$ \cite{AbMa78}.

If we substitute $\kappa $ and $\tau $ with their expressions in terms of
the variables $q,r$ and $s$, then $H$ in Eq.(\ref{HamLu1}) has the form
\begin{equation}
H\left( s,p,q,r\right) =e^{(\alpha +\gamma )t}\left[ (q+s)\arcsin \left(
\frac{q}{\sqrt{q^{2}+r^{2}}}\right) -(\alpha -1)r\right] +\frac{1}{2}%
p^{2}e^{-\alpha t},  \label{LuI3}
\end{equation}%
which is the Hamiltonian function of the nonautonoumous system (\ref{Lu3})
with a degenerate Poisson matrix having three-component vectors $\mathbf{U}%
=(0,-r,0)$ and $\mathbf{V}=(r,0,q)$. Since the L\"{u} system (\ref{LunonAut}%
) is the nonautonoumous the Hamiltonian $H$ in Eq.(\ref{LuI3}) is not an
integral invariant of the L\"{u} system.

To have a time-independent first integral of the L\"{u} system (\ref%
{LunonAut}), we choose $\gamma=-2\alpha$ and define a new time variable $%
\bar {t}=-e^{-at}/a$. This enables us to write the system in an autonomous
form%
\begin{equation}
p^{\prime}=\alpha
q+s,\;\;q^{\prime}=-rp\;\;\;r^{\prime}=qp\;\;\;s^{\prime}=rp
\end{equation}
where prime denotes derivative with respect to the time parameter $\bar{t}$.
In this case, the reduced dynamics in Eq.(\ref{Lu3}) becomes autonomous
\begin{equation}
p^{\prime}=\kappa+(\alpha-1)q,\;\;\;q^{\prime}=-p\sqrt{\tau-q^{2}}.
\label{Lu2aut}
\end{equation}
Integrating this we find that%
\begin{equation*}
H=\frac{1}{2}p^{2}+\kappa\arcsin\left( \frac{q}{\sqrt{\tau}}\right)
-(\alpha-1)\sqrt{\tau-q^{2}}
\end{equation*}
is the Hamiltonian of the system (\ref{Lu2aut}) and it is conserved.
Substituting $H_{1}=\kappa$ and $H_{2}=\tau$ from Eq.(\ref{ILu}), we arrive
at the third autonomous conserved quantity
\begin{equation}
H_{3}\left( s,p,q,r\right) =\frac{1}{2}p^{2}+\left( q+s\right) \arcsin\left(
\frac{q}{\sqrt{q^{2}+r^{2}}}\right) -(\alpha-1)r,  \label{H3Ray}
\end{equation}
of the L\"{u} system.

Note that, we started with L\"{u} system in Eqs.(\ref{LuG}) and showed that
when the parameters are satisfying the relations $\gamma=-\beta=\delta
=-2\alpha$, the system has three time-independent first integrals, namely $%
H_{1}$, $H_{2}$ in Eq.(\ref{ILu}) and $H_{3}$ in Eq.(\ref{H3Ray}). Now, we
follow the definitions in Eqs.(\ref{N}) in order to obtain the
three-component vector functions
\begin{align*}
\mathbf{U}_{1} & =(-2\alpha q+s+r\arcsin(\frac{q}{\sqrt{q^{2}+r^{2}}}%
),-rp,qp),\text{ }\mathbf{V}_{1}=-2\arcsin(\frac{q}{\sqrt{q^{2}+r^{2}}}%
)\left( 0,q,r\right) , \\
\mathbf{U}_{2} & =(\alpha-1+\left( q+s\right) \frac{q}{q^{2}+r^{2}},0,p),%
\text{ }\mathbf{V}_{2}=(-p,-\left( q+s\right) \frac{r}{q^{2}+r^{2}}%
,\alpha-1+\left( q+s\right) \frac{q}{q^{2}+r^{2}}), \\
\mathbf{U}_{3} & =2r(1,0,0),\text{ }\mathbf{V}_{3}=2(0,q,r),
\end{align*}
of mutually compatible Poisson structures $N^{\left( 1\right) }$, $N^{\left(
2\right) }$ and $N^{\left( 3\right) }$, respectively. The common conformal
factor is $-1/2$.

\subsection{Hyperchaotic Qi system \label{FIHQS}}

The hyperchaotic Qi system is given by the following set equations
\begin{align}
\dot{u} & =-\delta u+\lambda z+xy,  \notag \\
\dot{x} & =\alpha(y-x)+yz  \notag \\
\dot{y} & =\beta(x+y)-xz  \notag \\
\dot{z} & =-\gamma z-\epsilon u+xy  \label{QiSystem}
\end{align}
where $\alpha,\beta,\gamma,\epsilon,\delta$ and $\lambda$ are real constants
called the parameters. When the parameters satisfy the constraints
\begin{equation}
\alpha+\beta=0,\;\;\;\gamma+\epsilon+\lambda=\delta,  \label{paraQi}
\end{equation}
Qi system admits the following time-dependent first integrals
\begin{equation}
I_{1}=(z-u)e^{(\gamma+\lambda)t},\;\;\;I_{2}=(x^{2}+y^{2})e^{2\alpha t}.
\label{QiFI}
\end{equation}

We introduce a transformation $\left( u,x,y,z\right) \rightarrow\left(
s,p,q,r\right) $ given by%
\begin{equation}
s=ue^{(\gamma+\lambda)t},\;\;q=ye^{\alpha t},\text{ \ \ }p=xe^{\alpha
t},\;\;r=ze^{(\gamma+\lambda)t}.  \label{Qitrans}
\end{equation}
In the new coordinates, the system assumes the form
\begin{align}
\dot{s} & =\lambda r-\epsilon s+pqe^{(\gamma+\lambda-2\alpha)t}
\label{hyperQid} \\
\dot{q} & =p(\beta-re^{-(\gamma+\lambda)t})  \notag \\
\dot{p} & =q(re^{-(\gamma+\lambda)t}-\beta)  \notag \\
\dot{r} & =\lambda r-\epsilon s+pqe^{(\gamma+\lambda-2\alpha)t}  \notag
\end{align}
while the first integrals $I_{1}$ and $I_{2}$ in Eq.(\ref{QiFI}) becomes
autonomous
\begin{equation}
H_{1}=r-s,\;\;\;H_{2}=p^{2}+q^{2}.  \label{QiFI2}
\end{equation}
Note that, time derivatives $\dot{r}$ and $\dot{s}$ are the same, this is a
manifestation of functional form of the first integral $H_{1}=r-s$. On the
intersection of the level surfaces of $H_{1}=\kappa$ and $H_{2}=\tau$ the
non-autonomous system (\ref{hyperQid}) reduces to the planar system
\begin{equation}
\dot{r}=\epsilon\kappa+\left( \lambda-\epsilon\right) r+q\sqrt{\tau-q^{2}}%
e^{(\gamma+\lambda-2\alpha)t},\text{ \ \ }\dot{q}=\sqrt{\tau-q^{2}}%
(\beta-re^{-(\gamma+\lambda)t})  \label{QiPl1}
\end{equation}
upon elimination of the variables $p$ and $s$ in favor of $q$ and $r$. The
latter system of planar ODEs admits a Jacobi last multiplier given by
\begin{equation}
M=\frac{e^{(\epsilon-\lambda)t}}{\sqrt{\tau-q^{2}}}.  \label{hyperQiJLM}
\end{equation}
The multiplier $M$ is time-dependent, hence we introduce two auxiliary
functions%
\begin{equation}
\psi(q,r,t)=(\lambda-\epsilon)r,\text{ \ \ }\phi(q,r,t)=\beta\sqrt{\tau-q^{2}%
},  \label{QiAux1}
\end{equation}
satisfying the identity in Eq.(\ref{5.12}). To find Hamiltonian of the
reduced system in Eq.(\ref{QiPl1}), we recall the defining equation (\ref%
{5.13}), hence compute the Hamiltonian function as
\begin{equation}
H\left( q,r\right) =e^{(\epsilon-\lambda)t}\left( e^{-(\gamma+\lambda )t}%
\frac{r^{2}}{2}+\epsilon\kappa\arcsin\left( \frac{q}{\sqrt{\tau}}\right) +%
\frac{q^{2}}{2}e^{\left( \lambda+\gamma-2\alpha\right) t}\right) .
\label{H2DQi1}
\end{equation}
Note that, the system (\ref{QiPl1}) is nonautonomos, hence $H$ is not a
constant of motion. Using Eq.(\ref{Canoor}), we arrive at the canonical
coordinates
\begin{equation*}
Q=e^{(\epsilon-\lambda)t}r,\text{ \ \ }P=\arcsin\left( \frac{q}{\sqrt{\tau}}%
\right) -\beta t.
\end{equation*}
In the canonical coordinates, the planar system (\ref{QiPl1}) takes the form
\begin{equation}
\dot{Q}=e^{(\epsilon-\lambda)t}\epsilon\kappa+\frac{\tau}{2}\sin\left(
2\left( P+\beta t\right) \right) e^{\left( \epsilon+\gamma-2\alpha\right) t},%
\text{ \ \ }\dot{P}=-Qe^{-\left( \epsilon+\gamma\right) t},  \label{Qi2DCan}
\end{equation}
whereas the Hamiltonian for the planar system (\ref{H2DQi1}) assumes the
form
\begin{equation}
H\left( Q,P\right) =e^{-(\gamma+\epsilon)t}\frac{Q^{2}}{2}+\epsilon\kappa
e^{(\epsilon-\lambda)t}\left( P+\beta t\right) +\frac{\tau e^{\left(
\epsilon+\gamma-2\alpha\right) t}}{2}\sin^{2}\left( P+\beta t\right) .
\label{HamQi2DCan}
\end{equation}
So that, the canonical system (\ref{Qi2DCan}) can be recasted as $\dot {Q}%
=\partial H/\partial P$ and $\dot{P}=-\partial H/\partial Q$. Note that, the
Hamiltonian function in Eq.(\ref{HamQi2DCan}) is not conserved, that is $%
dH/dt=\partial H/\partial t$, because the system (\ref{Qi2DCan}) is
nonautonomous. If we substitute $\kappa$ and $\tau$ with their expressions
in terms of the variables $q,r$ and $s$, then the Hamiltonian in (\ref%
{H2DQi1}) becomes time-dependent function%
\begin{equation}
H\left( q,p,r,s\right) =e^{(\epsilon-\lambda)t}\left( e^{-(\delta +\lambda)t}%
\frac{r^{2}}{2}+\epsilon\left( r-s\right) \arcsin\left( \frac {q}{\sqrt{%
q^{2}+p^{2}}}\right) +\frac{q^{2}}{2}e^{\left( \lambda
+\gamma-2\alpha\right) t}\right) .  \label{QiH31}
\end{equation}

In order to have a time-independent integral of the Qi system we take $%
\lambda=-\gamma$ and $\alpha=0$ and change time variable $t$ with $\bar {t}=%
\frac{1}{\epsilon-\lambda}e^{\left( \epsilon-\lambda\right) t}$. In these
circumstances, the reduced planar ODE system in Eq.(\ref{QiPl1}) takes the
form%
\begin{equation}
r^{\prime}=\frac{1}{\bar{t}}\left( \frac{\epsilon\kappa}{\epsilon-\lambda }+q%
\sqrt{\tau-q^{2}}-r\right) ,\text{ \ \ }q^{\prime}=-\frac{r}{\bar{t}}\sqrt{%
\tau-q^{2}},
\end{equation}
where prime denotes derivative with respect to the time parameter $\bar{t}$.
\ In this case, the Jacobi's last multiplier becomes $M=\bar{t}/\sqrt {%
\tau-q^{2}}$ and the auxiliary functions turn out to be $\psi=-r/\bar{t}$
and $\phi=0$. Using Eq.(\ref{5.13}), we compute
\begin{equation*}
H\left( r,q\right) =\frac{1}{\epsilon-\lambda}\left( \epsilon\kappa
\arcsin\left( \frac{q}{\sqrt{\tau}}\right) +\frac{q^{2}}{2}+\frac{r^{2}}{2}%
\right) ,
\end{equation*}
which does not depend on time explicitly and hence after the substitution of
the first two integrals $\kappa=r-s$ and $\tau=q^{2}+p^{2}$ we arrive at the
new time-independent integral%
\begin{equation}
H_{3}\left( q,p,r,s\right) =\frac{1}{\epsilon-\lambda}\left( \epsilon \left(
r-s\right) \arcsin\left( \frac{q}{\sqrt{q^{2}+p^{2}}}\right) +\frac{q^{2}}{2}%
+\frac{r^{2}}{2}\right) .  \label{QiH32}
\end{equation}
of the Qi system.

We thus have obtained three time-independent integrals of the Qi system,
namely $H_{1}$, $H_{2}$ in Eq.(\ref{QiFI2}) and $H_{3}$ in Eq.(\ref{QiH32}),
with the parameters satisfying the conditions $\lambda =-\gamma ,$ $\delta
=\epsilon $ and $\alpha ,\beta =0.$ Using the definitions in Eq.(\ref{N}),
the tri-Hamiltonian structure of the Qi system is obtained as follows. For
the Hamiltonian function $H_{1}=r-s$, the three-component vectors of the
corresponding Poisson structure $N^{\left( 1\right) }$ are given by
\begin{align*}
\mathbf{U}_{1}& =\frac{2}{\epsilon -\lambda }(\epsilon p\arcsin (\frac{q}{%
\sqrt{q^{2}+p^{2}}})+pr),-\epsilon q\arcsin (\frac{q}{\sqrt{q^{2}+p^{2}}}%
)-qr,-\epsilon \left( r-s\right) -pq) \\
\mathbf{V}_{1}& =\frac{2\epsilon }{\epsilon -\lambda }\arcsin (\frac{q}{%
\sqrt{q^{2}+p^{2}}})(q,p,0).
\end{align*}%
For the Hamiltonian function $H_{2}=q^{2}+p^{2}$, the three-component
vectors of $N^{\left( 2\right) }$ are
\begin{align*}
\mathbf{U}_{2}& =\frac{1}{\epsilon -\lambda }(-\epsilon \left( r-s\right)
\frac{q}{p^{2}+q^{2}},-\epsilon \left( r-s\right) \frac{p}{q^{2}+p^{2}}-q,0)
\\
\mathbf{V}_{2}& =\frac{1}{\epsilon -\lambda }(\epsilon \left( r-s\right)
\frac{p}{q^{2}+p^{2}}+q,-\epsilon \left( r-s\right) \frac{q}{p^{2}+q^{2}},r).
\end{align*}%
Finally, for the Hamiltonian $H_{3}$ in Eq.(\ref{QiH32}), we have
\begin{equation*}
\mathbf{U}_{3}=\left( -2p,2q,0\right) ,\text{ \ \ }\mathbf{V}_{3}=\left(
-2q,-2p,0\right) ,
\end{equation*}%
which may be related to four dimensional Lie algebras \cite{LaPi12}.

\section{Discussion and outlook}

In this paper, after presenting some technical details about
superintegrability and Poisson structures in four dimensions, we have
studied superintegrable and tri-Hamiltonian structures of a set of four
dimensional dynamical systems, namely Shivamoggi and generalized
Raychaudhuri equations, hyperchaotic L\"{u} and Qi systems. Generically
hyperchaotic attractors are not expected to be integrable, but in this
paper, we have showed that, for some particular values of their parameters,
hyperchaotic L\"{u} and Qi systems exhibit not only integrable but
superintegrable properties admitting tri-Hamiltonian formulations.

\section*{Acknowledgement}

We gratefully acknowledge support from Professor G. Rangarajan and National
Mathematics Initiative programme at IISC Mathematics Department where the
work was started. HG thanks Yeditepe University for travel support.

\end{document}